\documentclass[12pt]{amsart}
\usepackage{amssymb, color, graphics, colordvi}
\usepackage{multirow}

\numberwithin{equation}{section}

\newtheorem{thm}{Theorem}
\numberwithin{thm}{section}
\newtheorem{prop}[thm]{Proposition}
\newtheorem{lemma}[thm]{Lemma}

\newtheorem{example}[thm]{Example}
\newtheorem{remark}[thm]{Remark}
\newtheorem{defi}[thm]{Definition}

\newcounter{FNC}[page]
\def\newfootnote#1{{\addtocounter{FNC}{2}$^\fnsymbol{FNC}$%
     \let\thefootnote\relax\footnotetext{$^\fnsymbol{FNC}$#1}}}

\pagespan{1}{60}
\copyrightinfo{}{}
\newcommand{\calA}{\mathcal{A}}

\newcommand{\calC}{\mathcal{C}}

\newcommand{\calS}{\mathcal{S}}
\newcommand{\calR}{\mathcal{R}}

\newcommand{\C}{\mathbb{C}}
\newcommand{\N}{\mathbb{N}}
\newcommand{\R}{\mathbb{R}}
\newcommand{\Z}{\mathbb{Z}}

\title{Polynomial systems supported on circuits
and dessins d'enfants}

\author{Frederic Bihan}
\address{Laboratoire de Math\'ematiques\\
         Universit\'e de Savoie\\
         73376 Le Bourget-du-Lac Cedex\\
         France}

\email{Frederic.Bihan@univ-savoie.fr}

\begin{document}

\begin{abstract}
We study polynomial systems whose equations
have as common support a set $\calC$ of $n+2$ points in ${\Z}^n$
called a circuit.
We find a bound on the number of real
solutions to such systems which depends on $n$, the dimension of the affine span of
the minimal affinely dependent subset of $\calC$,
and the rank modulo $2$ of $\calC$. We prove that this bound is sharp
by drawing so-called dessins d'enfant on the Riemann sphere.
We also obtain that the maximal number of solutions
with positive coordinates to systems supported on circuits in ${\Z}^n$ is $n+1$, which
is very small comparatively to the bound given by the Khovanskii fewnomial theorem.
\end{abstract}
\maketitle

\section*{Introduction and statement of the main results}

The support of a multivariate polynomial is the set of exponent vectors
of the monomials appearing in the polynomial.
A polynomial system of $n$ polynomial equations in $n$ variables is supported on $\calA \subset {\Z}^n$
if $\calA$ is the common support of each polynomial in the system.
A theorem due to Kouchnirenro gives that the number of (simple) solutions in the complex
torus $({\C}^*)^n$ to a generic polynomial system supported on $\calA$ is the volume $v(\calA)$
of the convex hull of $\calA$,
normalized so that the unit cube $[0,1]^n$ has volume $n!$.

Here, we consider generic polynomial systems with real coefficients, and are interested in their numbers
of real solutions in the real torus $({\R}^*)^n$.
Contrary to the complex case, the number of real solutions depends on the coefficients of the system, and one of the most important
question is to find a sharp upper bound (for related results, see \cite{BeBiS,LRW03, So, So1} for example). A trivial bound is given by $v(\calA)$,
the number of complex solutions.
Another bound due to Khovanskii depends only on the number $n$ of variables
and the cardinality of the support $\calA$. The above Kuchnirenko result shows immediately
that such a bound belongs to the ``real world'', that is, cannot exist as a bound on the number of complex solutions.
The Khovanskii bound is
$$2^n 2 ^{\binom{|\calA|}{2}}
 \cdot (n+1)^{|\calA|},$$
and gives an easy way to construct supports $\calA$ for which
the Kuchnirenko bound $v(\calA)$ is not sharp.
The Khovanskii bound can be stated alternatively without the term $2^n$
as a bound on the number of {\it positive solutions}, which are solutions with positive coordinates.

Since we are interested only in the solutions in the real torus, we have the freedom to translate the support $\calA$
by an integral vector and also to choose a basis of the lattice ${\Z}^n$.
A translation by an integral vector corresponds to multiplying each polynomial of the system by a monomial, while a change
of basis for the lattice ${\Z}^n$ corresponds to a monomial change of coordinates (with coefficients $1$) for the torus.
These operations do not change the number of positive, real and complex solutions. In particular, we will always assume that
$0 \in \calA$. The {\it rank modulo 2} of $\calA$ is the rank of the reduction modulo $2$
of the matrix obtained by putting in columns the non zero elements of $\calA$. We will denote it by
$rk \, (\bar{\calA})$.
%
\smallskip

If $\calA \subset {\Z}^n$ and $|\calA| <n+1$, then $v(\calA)=0$, so that supports with such numbers of elements
are not interesting. The first case where $v(\calA)$ is different from $0$ arises when $|\calA|=n+1$,
and the convex hull of $\calA$ is an $n$-dimensional simplex. This case is easy. The number of positive solutions
to such systems is at most $1$, so that the number of real solutions is at most $2^n$.
In fact, the maximal number of real solutions is $2^{n-rk \, (\bar{\calA})}$ (\cite{S}, see also Lemma~\ref{L:simplex}).
\smallskip

The goal of this paper is to obtain similar results in the first non trivial case, when the support is {\it a circuit}.
A circuit is a set $\calC \subset {\Z}^n$ of $n+2$ points which affinely span ${\R}^n$.
It contains an unique minimal affinely dependent subset, and we denote by $m(\calC)$ the dimension
of the affine span of this subset. If $m(\calC)=n$, that is, if any proper subset of $\calC$ is affinely independent,
then $\calC$ is called a {\it non degenerate circuit}. We have $m(\calC)=1$ for example
when three points of $\calC$ lie on a same line.
In~\cite{BeBiS}, the authors consider circuits $\calC$ (and also more general supports called near circuits)
such that the index of ${\Z}\calC$ in ${\Z}^n$ is odd. This corresponds to circuits $\calC \subset {\Z}^n$
with $rk \, (\bar{\calC})=n$. They obtain the upper bound $2m(\calC)+1$ on the number of real solutions
and also prove that this bound is sharp among systems supported on circuits $\calC \subset {\Z}^n$
with $rk \, (\bar{\calC})=n$. In particular, this gives the sharp bound $2n+1$ which can be attained
only with non degenerate circuits. We generalize the results in~\cite{BeBiS}
to arbitrary values of the rank modulo $2$. The first result gives sharp bounds on the number of positive solutions.
Obviously, multiplying by $2^n$ such a bound gives a bound on the number of real solutions,
which is sharp among systems supported on circuits $\calC$ with $rk \, (\bar{\calC})=0$.
\smallskip

\noindent {\bf Theorem A.}
{\it The number of positive solutions to a generic real polynomial system supported on a
circuit $\calC \subset {\Z}^n$ is at most
$$m(\calC)+1.$$
Therefore, the number of real solutions to a generic real polynomial system supported on a
circuit $\calC \subset {\Z}^n$ is at most
$$2^n(m(\calC)+1).$$
Moreover, the first bound, and thus the second bound, is sharp.
Namely, for any integer $m$ such that $1 \leq m \leq n$, there exist a circuit $\calC \subset {\Z}^n$ with $m(\calC)=m$,
and a system supported on $\calC$ which has $m+1$ positive solutions.

As a consequence, the number of positive solutions to a generic system supported on a circuit in ${\Z}^n$ is at most
$$n+1$$
while its number of real solutions is at most
$$2^n \cdot (n+1)$$
and these bounds are sharp and can be attained only with non degenerate circuits.}
\smallskip

\noindent In particular, as in the case of supports forming a simplex, the Khovansky bound is far from being sharp
among systems supported on circuits. Theorem A follows from Proposition~\ref{P:boundpositive} and Theorem~\ref{T:positivesharp}.
\smallskip

\noindent {\bf Theorem B.}
{\it The number, $N$, of real solutions to a generic real polynomial system supported
on a circuit $\calC \subset {\Z}^n$ satisfies
\begin{enumerate}
\item If $ {{rk \, ({\bar \calC})}} \leq m(\calC)$, then
$$N \leq 2^{n-{{rk \, ({\bar \calC})}}} \cdot \left(m(\calC)+ {{rk \, ({\bar \calC})}}+1\right).$$
\item if $ {{rk \, ({\bar \calC})}} \geq m(\calC)$, then
$$N \leq 2^{n-{{rk \, ({\bar \calC})}}} \cdot \left(2m(\calC)+1\right).$$
\end{enumerate}
Moreover, both these bounds are sharp. Namely, let $n,m,R$ be integers such that $1 \leq m \leq n$
and $0 \leq R \leq n$. If $R \leq m$ (resp. $R \geq m$), there exist a circuit
$\calC \subset {\Z}^n$ with $m(\calC)=m$, ${{rk \, ({\bar \calC})}}=R$,
and a system supported on $\calC$  whose number of real solutions is the bound in  
(1) (resp. the bound in (2)).}
\smallskip

The bounds in Theorem B look like the bound $2^{n-rk \, (\bar{\calA})}$ when $\calA$ forms an
$n$-dimensional simplex.
However, there is an essential difference in that the sharp bounds in Theorem B do not provide
the maximal number of real solutions to systems supported on a {\it given} circuit.
Theorem B follows from Theorem~\ref{T:upper bound} and Theorem~\ref{T:sharp}.
\bigskip

We consider the eliminant defined in~\cite{BeBiS} of a system supported on a circuit $\calC \subset {\Z}^n$.
Assuming that $\calC=\{0,\ell e_1,w_1,\dotsc,w_n\} \subset {\Z}^n$, this eliminant is a univariate polynomial
of the form
$$f(x)=x^{\lambda_0}\prod_{i=1}^t (g_i(x))^{\lambda_i} - \prod_{i=t+1}^{\nu} (g_i(x))^{\lambda_i},$$
where each polynomial $g_i$ has the form $g_i(x)=a_i+b_ix^{\ell}$ with $a_i$ and $b_i$ non zero real numbers.
Here, the coefficients $a_i$ and $b_i$ come from the coefficients of the system while the other numbers are determined
by $\calC$. The number $\nu$ is equal to $m(\calC)$ if $0$ and $\ell e_1$ belong to the minimal
affinely dependent subset of $\calC$. If $rk \, (\bar{\calC}) \neq 0$, then the integer $\ell$ can be assumed to be odd.
We choose to distinguish the case $rk \, (\bar{\calC}) \neq 0$ from the case $rk \, (\bar{\calC})=0$
for which the number of real solutions is given by the number of positive ones.
When $rk \, (\bar{\calC})=n$, it is proved in~\cite{BeBiS} that the real solutions to the system
are in bijection with the real roots of $f$ via the projection onto the first coordinate axis.
In general, the real solutions project onto the real roots of $f$ which satisfy sign conditions
involving products of polynomials $g_i$. Moreover, such a real root is the image of $2^{n-{{rk \, ({\bar \calC})}}}$
solutions to the system, at least when $\ell$ is odd. The positive solutions project bijectively
onto the positive roots of $f$ at which all $g_i$ are positive.
We determine the above sign conditions, and prove the upper bounds in Theorem A and B using essentially Rolle's theorem.
\smallskip

Writing the eliminant as $f=P-Q$, we see that the number of real (and positive) solutions is closely related to
the arrangement of the roots of $P$, $Q$ and $f$,  their multiplicities being determined by $\calC$.
We consider the rational function $$\phi=f/Q=P/Q-1: {\C}P^1 \rightarrow {\C}P^1.$$
The roots of $P$, $Q$ and $f$ are the inverse images of $-1$, $\infty$ and $0$, respectively.
These roots lie on the graph 
$$\Gamma=\phi^{-1}({\R}P^1) \subset {\C}P^1$$
and we see $\Gamma$ as an example of so-called {\it dessin d'enfant}.
We use then the observation made in~\cite{Br, O} (see also~\cite{NSV}) that polynomials $P$, $Q$ and $f=P-Q$ with prescribed
arrangement and multiplicities of their real roots can be constructed drawing so-called {\it real rational
graphs} on ${\C}P^1$. The sharpness of the bounds in Theorem A and B is then proved by means of such graphs.

The paper is organized as follows.
In Section 1, we define the eliminant $f$ of a system
and explain the relation between its real roots and the real (resp. positive) solutions to the system.
In Section 2, we prove the upper bounds in Theorem A and B, while in Section 3
we introduce real rational graphs and achieve constructions proving the sharpness of these upper bounds.
\bigskip

\noindent
{\bf Acknowledgements.} I'm greatly indebted to Ilia Itenberg and Frank Sottile for their support
during these last years. Thanks to Frank Sottile and Erwan Brugall\'e for useful comments
on the first version of this paper. Thanks also to Benoit Bertrand and J. Maurice Rojas
for their interest in this work.
\bigskip

\section{Elimination}\label{S:elimination}

A real polynomial system supported on a circuit
$$\calC=\{w_{-1},w_0,\cdots,w_n\} \subset {\Z}^n$$
is called {\it generic} if it has $v(\calC)$ solutions in ${({\C}^*)}^n$. This forces each solution
to be a simple solution. We are interested in the real solutions to such a system, by which we mean
solutions in the real torus ${({\R}^*)}^n$.

Translating $\calC$ by an integral vector and choosing a basis for the lattice ${\Z}^n$ if necessary,
we may assume without loss of generality that $w_{-1}=0$, and $w_0=\ell e_1$ for some
positive integer $\ell$. Perturbing slightly the system, and using Gaussian elimination,
it becomes equivalent to
a system of the form
\begin{equation}\label{E:reduced}
S \; : \qquad x^{w_i}=g_i(x_1) \;, \quad i=1,\dotsc,n,
\end{equation}
where $g_i(x_1)=a_i+b_ix_1^{\ell}$ with $a_i,b_i$ non zero real numbers for $i=1,\dotsc,n$ and
$g_1,\dotsc,g_n$ have distinct roots.
Reordering the vectors $w_1,\dotsc,w_n$ if necessary,
the affine primitive relation on $\{0,e_1,w_1,\dotsc, w_n\}$ can be written as
\begin{equation}\label{E:primitiveaffinerelation}
\lambda_0 \, e_1+\sum_{i=1}^t \lambda_iw_i = \sum_{i=t+1}^{\nu} \lambda_i w_i,
\end{equation}
where $1 \leq \nu \leq n$, $\lambda_0,\dotsc,\lambda_{\nu}$ are coprime integers, $\lambda_0 \geq 0$
and $\lambda_1,\dotsc,\lambda_{\nu} >0$. Note that here we could have $t=0$ or $t=\nu$
so that one of the two sums collapses to $0$.
The integer
$$\delta:= \lambda_0+\sum_{i=1}^t \lambda_i-\sum_{i=t+1}^{\nu} \lambda_i$$
is (up to sign) the coefficient of $w_{-1}=0$ in ~\eqref{E:primitiveaffinerelation}
(the term $\delta w_{-1}$ does not appear in~\eqref{E:primitiveaffinerelation} since it is $0$).
Multiplying~\eqref{E:primitiveaffinerelation} by $\ell$ gives an affine relation on $\calC$.
Set
\begin{equation}\label{E:minimal dimension} 
m(\calC):=\nu -\chi(\lambda_0=0)-\chi(\delta=0)
\end{equation}
where $\chi(Y)$ is the boolean truth value of $Y$: $\chi(Y)=1$ if $Y$ is true and $\chi(Y)=0$ otherwise.
Then $m(\calC)+2$ is the number of non zero coefficients in the primitive affine relation on $\calC$, so that
$m(\calC)$ is the dimension of the affine span of the minimal affinely dependent subset of $\calC$.
The case $m(\calC)=n$ arises when $\calC$ is a non degenerate circuit, that is,
when any proper subset of $\calC$ is affinely independent.

Define the {\it eliminant} of $S$ to be the following univariate polynomial
\begin{equation}\label{E:eliminant}
f(x_1)=x_1^{\lambda_0}\prod_{i=1}^t (g_i(x_1))^{\lambda_i} - \prod_{i=t+1}^{\nu} (g_i(x_1))^{\lambda_i}.
\end{equation}

From~\eqref{E:reduced} and~\eqref{E:primitiveaffinerelation}, it follows immediately
that if $x=(x_1,\dotsc,x_n)$ is a solution to $S$ then $x_1$ is a root of the eliminant $f$.

Write each vector $w_i=l_ie_1+v_i$ with $v_i \in {\Z}^{n-1}$. We get the relations
\begin{equation}\label{E:primitiveaffinerelationv}
\sum_{i=1}^t \lambda_iv_i = \sum_{i=t+1}^{\nu} \lambda_i v_i
\end{equation}
and
\begin{equation}\label{E:primitiveaffinerelationl}
\lambda_0+\sum_{i=1}^t \lambda_il_i = \sum_{i=t+1}^{\nu}\lambda_i l_i.
\end{equation}

The fact that $\lambda_0,\dotsc,\lambda_{\nu}$ are coprime and~\eqref{E:primitiveaffinerelationl}
gives immediately that $\lambda_1,\dotsc,\lambda_{\nu}$
are also coprime, so that the relation~\eqref{E:primitiveaffinerelationv} is in fact the primitive affine relation
on $\{0,v_1,\dotsc,v_n\} \subset {\Z}^{n-1}$. In particular, at least one integer among $\lambda_1,\dotsc,\lambda_{\nu}$
should be odd.

\begin{prop}[\cite{BeBiS}]\label{P:Systemreduced}
Assume that $\lambda_j$ is odd with $j \in \{1,\dotsc,{\nu}\}$.
Then $x=(x_1,\dotsc,x_n)$ is a real solution to the system $S$ if and only if $x_1$ is a real root
of $f$ and $(x_2,\dotsc,x_{n})$ is a real solution to the system
$$y^{v_i}=g_i(x_1)/x_1^{l_i} \;, \quad i \in \{1,\dotsc,n\} \setminus \{j\}.$$
\end{prop}

For simplicity, we will often assume that $\lambda_1$ is odd and thus apply Proposition~\ref{P:Systemreduced}
with $j=1$. Define 
$$B:=(v_1,
v_2,  \dotsc , v_{n}),$$
the $n-1$ by $n$ matrix whose columns are the vectors $v_1,v_2,\dotsc,v_{n}$.
Let $\bar{B} \in M_{n-1,n}(\Z/2)$ be the
the reduction modulo $2$ of $B$ and $rk \,(\bar{B})$ be the rank of $\bar B$.

\begin{prop}\label{P:real}
The number of real solutions to $S$ is equal to $2^{n-1-rk \,(\bar{B})}$ times
the number of real roots $r$ of $f$ subjected to the sign conditions
$${(g_1(r)/r^{l_1})}^{\epsilon_1} \cdots {(g_{n}(r)/r^{l_{n}})}^{\epsilon_{n}}>0
$$
for any $\epsilon=(\epsilon_1,\dotsc,\epsilon_n)$ such that $\bar{\epsilon} \in \mbox{Ker} \, \bar{B}$.
\end{prop}
\noindent It is worth noting that the number of sign conditions can be reduced to $n-rk \,(\bar{B})$
using a basis of $\mbox{Ker} \, \bar{B}$.
\smallskip

Before giving the proof, we present the following well-known lemma (see \cite{S}).

\begin{lemma}\label{L:simplex}
Suppose that $0,u_1,\dotsc,u_n \in {\Z}^n$ form
an $n$-dimensional simplex , and let $\bar{U} \in M_n(\Z/2)$ denote the matrix
whose columns contains in this order the reductions modulo $2$ of $u_1,\dotsc,u_n$.
Then the number of real solutions in the torus to the system
\begin{equation}\label{E:intilasyst1}
x^{u_i}=c_i,\;, \quad i=1,\dotsc,n,
\end{equation}
where $c_i$ is a non zero real number,
is
\begin{enumerate}
\item $0$ or $2^{n-rk(\bar{U})}$ if the volume of the simplex is even,
\item $1$ if this volume is odd.
\end{enumerate}
Moreover, we have $2^{n-rk(\bar{U})}$ real solutions in the first case if and only if
$$c_1^{\epsilon_1}\cdots c_n^{\epsilon_n} >0$$
for any $(\epsilon_1,\dotsc,\epsilon_n) \in {\Z}^n$ whose reduction modulo $2$ belong to the kernel
of $\bar{U}$.
\end{lemma}
Note that the case $2)$ is in fact redundant since if the volume of the simplex is odd then
the kernel of $\bar{U}$ is reduced to $0$, hence $2^{n-rk(\bar{U})}=1$, and the sign conditions are empty.
\smallskip

\noindent {\bf Proof.}
There exist two matrices $L,R \in GL_{n}(\Z)$ and integers $a_1,\dotsc,a_{n}$
such that $LUR$ is the diagonal matrix $D=\mbox{diag} \, (a_1,\dotsc,a_{n})$.
Let $\tilde{u}_1,\dotsc,\tilde{u}_{n}$ be the columns vectors, in this order, of the matrix $LU$.
The matrix $L$ provides the multiplicative change of coordinates
$\tilde{y_i}= \prod_{j=1}^{n} y_j^{(L^{-1})_{ji}}$, $i=1,\dotsc,n$,
of the complex torus. This change of coordinates transforms our system to the system
$\tilde{y}^{\tilde{u}_i}=c_i  \;, \quad i=1,\dotsc,n$, which has the same number of real (and complex) solutions.
Multiplication on the right by $R$ transforms this system to the system

\begin{equation}\label{E:diagsyst}
\tilde{y_i}^{a_i}=d_i \quad \mbox{with}  \quad d_i:=\prod_{j=1}^{n}c_j^{R_{ji}}\;, \quad i=1,\dotsc,n.
\end{equation}
This system is equivalent to the previous one
since $\tilde{y_i}^{a_i}=\prod_{j=1}^{n}(\tilde{y}^{\tilde{u}_j})^{R_{ji}}$
and $R \in GL_{n}(\Z)$.
In particular, the initial system~\eqref{E:intilasyst1} and the system~\eqref{E:diagsyst}
have the same number of real solutions. This number of real solutions is
$0$ if $d_i<0$ for some even $a_i$, or $2^{n-rk(\bar{D})}$ if $d_i>0$ for all even $a_i$.
The conclusion is now obvious noting that $rk(\bar{D})=rk(\bar{U})$
and that the reductions modulo $2$ of the vector columns $(R_{1i},\dotsc,R_{n,i})$ for which $a_i$ is even
generate the kernel of $\bar{U}$. $\Box$
\smallskip

\noindent {\bf Proof of Proposition~\ref{P:real}}
Assume that $\lambda_1$ odd. According to Proposition~\ref{P:Systemreduced}, we have to count, for any root $r$ of $f$,
the number of real solutions to
\begin{equation}\label{E:intilasyst}
y^{v_i}=g_i(r)/r^{l_i} \;, \quad i=2,\dotsc,n.
\end{equation}
Set $c_i(r):=g_i(r)/r^{l_i}$, $i=1,\dotsc,n$.
Note that $0,v_2,\dotsc,v_n$ form an $(n-1)$-dimensional simplex for otherwise
there would be an additional affine relation on $\calC$ (so that the convex hull of $\calC$
would have dimension $<n$). Let $\bar{U} \in M_{n-1}(\Z/2)$ denote the matrix
whose columns contains in this order the reductions modulo $2$ of $v_2,\dotsc,v_n$.
According to Lemma~\ref{L:simplex}, the number of solutions to~\eqref{E:intilasyst} is either $0$, or $2^{n-1-rk(\bar{U})}$.
Moreover, it is $2^{n-1-rk(\bar{U})}$ if and only if
$$c_2(r)^{\epsilon_2} \cdots  c_n(r)^{\epsilon_n}>0$$
for any $(\epsilon_2,\dotsc,\epsilon_n)$ whose reduction modulo $2$ belongs to the kernel of
$\bar{U}=(\bar{v_2}, \; \dotsc, \bar{v_n})$.

It follows from the relation~\eqref{E:primitiveaffinerelationv}
that the reduction modulo $2$ of $\lambda:=(\lambda_1,\dotsc,\lambda_{\nu},
0,\dotsc,0)$ belongs to the kernel of $\bar B$. Moreover, since $\lambda_1$ is odd,
the kernel of $\bar B$ is the direct sum of $(\Z/2) \cdot \bar{\lambda}$ and the kernel (in $\{0\} \times (\Z/2)^{n-1}$) of
$\bar{U}$. This also implies that $rk(\bar{U})=rk(\bar{B})$. To finish, it remains to see that
$f(r)=0 \Rightarrow c_1(r)^{\lambda_1} \cdots c_{\nu}(r)^{\lambda_{\nu}} >0$.
$\Box$
\smallskip

\begin{remark}\label{R:redundant}
As explained at the end of the proof, the sign condition $c_1(r)^{\lambda_1} \cdots c_{\nu}(r)^{\lambda_{\nu}} >0$
given by the affine relation on $\{0,v_1,\dotsc,v_n\}$ is in fact redundant in Proposition~\ref{P:real}.
\end{remark}
\smallskip

A {\it positive} solution is a solution with positive coordinates.

\begin{prop}\label{P:positive}
The number of positive solutions to $S$ is equal to
the number of roots $r$ of $f$ subjected to the sign conditions
$$r>0 \; ,  \; g_1(r)>0 \; ,  \,\dotsc \, ,  \; g_n(r) >0.$$
\end{prop}
\smallskip

\noindent {\bf Proof.}
Assume that $\lambda_1$ is odd.
According to Proposition~\ref{P:Systemreduced} $(x_1,\dotsc,x_n)$ is a positive solution to $S$
if and only if $x_1$ is a positive root $r$ of $f$ and $(x_2,\dotsc,x_n)$ is a positive solution
to
\begin{equation}
y^{v_i}=g_i(r)/r^{l_i} \;, \quad i=2,\dotsc,n.
\end{equation}
As it is well-known, the previous system has a positive solution if and only if
$g_i(r)/r^{l_i}>0$ for $i=2,\dotsc,n$, and in this case it has exactly one positive solution.
To finish, it remains to note that $f(r)=0$, $r>0$ and $g_i(r)>0$ for $i=2,\dotsc,n$ implies that $g_1(r)>0$.
$\Box$
\smallskip

\section{Upper bounds}

Let $S$ be a system as in Section~\ref{S:elimination}
and write the eliminant~\eqref{E:eliminant} as
$$f=P-Q,$$
where $P(x_1)=x_1^{\lambda_0}\prod_{i=1}^t (g_i(x_1))^{\lambda_i}$ and
$Q(x_1)=\prod_{i=t+1}^{\nu} (g_i(x_1))^{\lambda_i}$.
Consider the rational function
\begin{equation}\label{E:rationalfunction}
\phi:=f/Q=P/Q-1.
\end{equation}
Recall that $m(\calC)$ is the dimension of the affine span of the minimal affinely dependent subset of $\calC$.

\begin{prop}\label{P:boundpositive}
The number of positive solutions to $S$ is no more than $m(\calC)+1$.
\end{prop}

\noindent {\bf Proof.}
Let ${\calR}_+$ denote the set of positive roots $r$ of $f$ satisfying $g_1(r)>0 \, ,  \,\dotsc \, ,  \,g_n(r) >0$
and let $N$ denote the number of elements of ${\calR}_+$.
The number of positive solutions to $S$ is equal to $N$ by Proposition~\ref{P:positive}.
Consider an open interval $I$ formed by two consecutive elements of
${\calR}_+$. This interval contains no root of $P$ and $Q$ (since each $g_i$ has simple real roots).
As a consequence, the function $\phi$ is bounded on $I$.
Moreover, $\phi$ takes the same value $0$ on the endpoints of $I$.
Therefore, by Rolle's theorem, the  derivative  $\phi'$ has at least one real root in $I$. 
This gives $N-1$ distinct positive roots of $\phi'$,
which are different from the roots of $P$ and $Q$. The roots of $\phi'$ are the roots of $P'Q-PQ'$. 
As noticed in (\cite{BeBiS}, Proof of Proposition 4.3.), the quotient

\begin{equation}\label{E:quotient}
\frac{P'(x)Q(x)-P(x)Q'(x)}{x^{\lambda_0-1}\prod_{i=1}^{\mu} (g_i(x))^{\lambda_i-1}},
\end{equation}
or, if $\lambda_0=0$, the same quotient but with  $\ell-1$ in place of $\lambda_0 -1$,
is a polynomial $H$ which can be written as $H(x)=h(x^{\ell})$ where $h$ is real polynomial of degree
$m(\calC)$ and with non zero constant term.
Thus the $N-1$ positive roots of $\phi'$ different from the roots of $P$ and $Q$ should be roots of $H$.
But $H$ has at most $m(\calC)$ distinct positive roots, hence $N-1 \leq m(\calC)$.
$\Box$
\smallskip

J. Maurice Rojas informed us that Proposition~\ref{P:boundpositive} can be obtained from~\cite{LRW03} (Lemma 2, Section 3).

Define the matrix $A$ as the matrix obtained
by putting (in this order) the vectors $w_{0},\dotsc,w_n$ in columns

\begin{equation}\label{E:matrixcircuit}
A:=(w_0,  \dotsc , w_{n})
 \in M_{n,n+1}(\Z)
\end{equation}
We have by definition ${{rk \, ({\bar \calC})}}={{rk \, ({\bar A})}}$, where $\bar{A}$ is the reduction modulo $2$ of $A$.
Recall that $w_{0}={\ell}e_1$ and $w_i=l_ie_1+v_i$ for $i=1,\dotsc,n$, so that

\begin{equation}\label{E:matrixcircuit2}
A=
\left(
\begin{array}{cccc}
{\ell} & l_1 & \dotsc & l_n \\
0 &     &       &      \\  
\vdots & v_1 & \dotsc & v_n \\
0 &     &       &     
\end{array}
\right) 
\end{equation}
The lower-right matrix $(v_1, \dotsc, v_n)$ is the matrix $B$ that we already defined.

\begin{thm}\label{T:upper bound}
Let $N$ be the number of real solutions to $S$.

\begin{enumerate}
\item If $ {{rk \, ({\bar \calC})}} \leq m(\calC)$, then
$$N \leq 2^{n-{{rk \, ({\bar \calC})}}} \cdot \left(m(\calC)+ {{rk \, ({\bar \calC})}}+1\right).$$
\item if $ {{rk \, ({\bar \calC})}} \geq m(\calC)$, then
$$N \leq 2^{n-{{rk \, ({\bar \calC})}}} \cdot \left(2m(\calC)+1\right).$$
\end{enumerate}
\end{thm}

\noindent {\bf Proof.}
Consider first the case ${{rk \, ({\bar \calC})}}=0$.
This corresponds to the case $\calC \subset 2{\Z}^n$.
Setting $\tilde{w}_i=w_i/2$, a system supported on $\calC$ can be rewritten as a system supported
on the circuit $\tilde{\calC}=\{0,\tilde{w}_{-1},\tilde{w}_0,\dotsc,\tilde{w}_n\}$ so that the number of
real solutions (in the real torus) to the system supported on $\calC$
is equal to $2^n$ times the number of positive solutions to the system supported on $\tilde{\calC}$.
Obviously, we have $m(\tilde{\calC})=m(\calC)$.
The bound $2^n(m(\calC)+1)$ for the number of real solutions to $S$
follows then from Proposition~\ref{P:boundpositive}.

Suppose now that ${{rk \, ({\bar \calC})}} \neq 0$, that is, $\calC$ is not contained in $2{\Z}^n$. 
Then, it is easy to show the existence of points $w_i,w_j \in \calC$ such that $w_i-w_j \notin 2{\Z}^n$, and $w_i$
belongs to the minimal affinely dependent subset of $\calC$.
Thus, reordering the points of $\calC$ so that $w_i$ becomes $w_{-1}$ and $w_j$ becomes $w_0$ if necessary,
we can assume that the number ${\ell}$ is an odd number, and the coefficient $\delta$ is non zero.
The fact that $\ell$ is odd gives that $\bar{A}$ and $\bar{B}$ have isomorphic kernels,
so that ${{rk \, ({\bar \calC})}}={{rk \, ({\bar A})}}=rk \,(\bar{B})+1$.

Define
$$c(r):=(c_1(r),\dotsc,c_n(r)) \quad \mbox{with} \quad
c_i(r):={g_{i}(r)/r^{l_{i}}} \;, \quad i=1,\dotsc,n.$$
Let $\pi_1,\dotsc,\pi_{n-{{rk \, ({\bar \calC})}}+1}$
be vectors in ${\Z}^{n-1}$ whose reductions modulo $2$ form
a basis of $\mbox{Ker} \, \bar{B}$. Then Proposition~\ref{P:real} gives
that the number of real solutions to $S$
is $2^{n-{{rk \, ({\bar \calC})}}}$  times the number of real roots $r$ of $f$ subjected to the sign conditions

\begin{equation}\label{E:cond1}
(c(r))^{\pi_i}>0 \quad \mbox{for} \quad i=1,\dotsc,n-rk \, ({\bar \calC})+1.
\end{equation}

Consider the ${n-{{rk \, ({\bar \calC})}}}+1$ by $n$ matrix $\Pi$ whose rows are given by the exponent vectors of the sign
conditions~\eqref{E:cond1}. Taking the Hermite normal form of the reduction modulo $2$ of $\Pi$, we obtain another
basis of $\mbox{Ker} \, \bar{B}$ producing equivalent sign conditions

\begin{equation}\label{E:normalcond}
\begin{array}{cccc}
c_{s_1}(r) \cdot \prod_{i>s_1, i \notin {\calS}} \, (c_i(r))^* & > & 0  \\
  & & &  \\
c_{s_2}(r)  \cdot \prod_{i>s_2, i \notin {\calS}} \, (c_i(r))^* & > & 0  \\
  & & &  \\ 
 \vdots & \vdots & \vdots \\ 
 & & & \\
c_{s_{n-{{rk \, ({\bar \calC})}}+1}}(r) \cdot \prod_{i>s_{n-{{rk \, ({\bar \calC})}}+1}, i \notin {\calS}} \, (c_i(r))^* & > & 0 
\end{array}
\end{equation}
where $1 \leq s_1 < s_2 < \dotsc < s_{n-{{rk \, ({\bar \calC})}}+1} \leq n$, ${\calS}=\{s_1,\dotsc,s_{n-{{rk \, ({\bar \calC})}}+1}\}$ and the stars
are exponents taking values $0$ or $1$.

It remains to show that the set ${\calR}$ of real roots $r$ of $f$
satisfying the sign conditions~\eqref{E:normalcond} has at most $m(\calC)+1+{{rk \, ({\bar \calC})}}$ elements if ${{rk \, ({\bar \calC})}} \leq m(\calC)$,
and at most $2m(\calC)+1$ elements if ${{rk \, ({\bar \calC})}} \geq m(\calC)$.
\smallskip

First note that the integer ${\ell}$ being odd, each $g_i$ has only one real root, that we denote by $\rho_i$. Moreover,
the polynomial $H$ defined in~\eqref{E:quotient} has at most $m(\calC)$ real roots since $H(x_1)=h(x_1^{\ell})$ with $h$ polynomial of degree $m(\calC)$.
Let $I$ denote an open (bounded) interval formed by two consecutive elements of ${\calR}$.
We claim that $I$ should contain at least one element among $0$,
$\rho_1,\dotsc,\rho_{\nu}$ and the real roots of $H$.
Indeed, if the function $\phi=f/Q$ is bounded on $I$, then $I$ should contain a real critical point (with finite critical value)
of $\phi$ by Rolle's theorem. The critical points of $\phi$ with finite critical values are contained in the set formed of $0$,
the roots $\rho_1,\dotsc,\rho_t$ of $P$ and the roots of $H$ (see the proof of Proposition~\ref{P:boundpositive}).
If $\phi$ is not bounded on $I$, then obviously $I$ should contain a root of $Q$, that is, a root $\rho_i$ for some
$i=t+1,\dotsc,\nu$.

There is at most one interval $I$ containing $0$
and at most $m(\calC)$ intervals $I$ containing a real root of $H$.
Let us concentrate now on the intervals $I$ which do no contain $0$ but contain a root
$\rho_i$ with $i=1,\dotsc,\nu$. The number of such intervals is obviously no more than $\nu$.
We claim that this number of intervals is also no more than ${{rk \, ({\bar \calC})}}-1$. Indeed, such an
interval $I$ should contain a $\rho_i$ with
$i \in \{1,\dotsc,n\} \setminus {\calS}$, for otherwise some sign condition
in~\eqref{E:normalcond} would be violated at the endpoints of $I$. The claim follows then as the cardinality
of $\{1,\dotsc,n\} \setminus {\calS}$ is ${{rk \, ({\bar \calC})}}-1$.

Therefore, the number of intervals $I$ formed by two consecutive elements of ${\calR}$ is no more than
$1+m(\calC)+{{rk \, ({\bar \calC})}}-1=m(\calC)+{{rk \, ({\bar \calC})}}$, and no more than $1+m(\calC)+\nu$. Obviously, the number of elements
of ${\calR}$ is equal to the number of these intervals $I$ added by one. The bound $m(\calC)+{{rk \, ({\bar \calC})}}$ on the number of
intervals $I$ gives thus
\begin{equation}\label{E:firstbound}
N \leq 2^{n-{{rk \, ({\bar \calC})}}} \cdot \left(m(\calC)+ {{rk \, ({\bar \calC})}}+1\right).
\end{equation}

Consider now the bound $1+m(\calC)+\nu$ on the number of intervals $I$. 
We see that if this bound is attained, then the polynomial $H$ has $m(\calC)$ real roots, and
each interval $I$ contains exactly one element among the real roots of $H$, the roots $\rho_1,\dotsc,\rho_{\nu}$,
and $0$. This forces each of these points to be a critical point of $\phi$ with even multiplicity (and
unbounded critical values for the roots of $Q$). The $m(\calC)$ real roots of $H$ are simple roots of $H$,
hence critical points with multiplicity $2$ of $\phi$. The multiplicity of $0$ with respect to $\phi $ is $\lambda_0$,
while the multiplicity of $\rho_i$, $i=1,\dotsc,\nu$, with respect to $\phi$ is $\pm \lambda_i$.
Consequently, if the number of intervals $I$ is $1+m(\calC)+\nu$, then $\lambda_0,\lambda_1,\dotsc,\lambda_{\nu}$ are even:
a contradiction since these numbers are coprime. Hence, the number of intervals $I$
is no more than $m(\calC)+\nu$. If $\lambda_0$ and $\delta$ are both non zero, then $\nu=m(\calC)$ so that
the number of intervals $I$ is no more than $2m(\calC)$. This gives the bound
\begin{equation}\label{E:secondbound}
N \leq 2^{n-{{rk \, ({\bar \calC})}}} \cdot \left(2m(\calC)+1\right).
\end{equation}
Recall that we have assumed at the beginning that $\delta \neq 0$.
Hence, it remains to see what happens when $\lambda_0=0$. In this case, we have
$m(\calC)=\nu-1$ and $0$ cannot be a critical point of $\phi$.
Arguing as before, we obtain then that the number of intervals $I$ is no more than $m(\calC)+\nu$.
Moreover, if this number is equal to $m(\calC)+\nu$, then $\lambda_1,\dotsc,\lambda_{\nu}$ are even:
a contradiction since $\lambda_1,\dotsc,\lambda_{\nu}$ are coprime. Hence, the number of intervals
$I$ is no more than $m(\calC)+\nu-1=2m(\calC)$, and we retrieve the bound~\eqref{E:secondbound}.
$\Box$
\smallskip

\section{Real rational graphs and sharpness of bounds}

\subsection{Real rational graphs}

This subsection comes from \cite{Br} (see also~\cite{Brt,O,NSV}).
Consider the  function~\eqref{E:rationalfunction}
$$\phi =f/Q=P/Q-1$$
as a rational function ${\C}P^1 \rightarrow {\C}P^1$. The degree of $\phi$ coincides with that of $f$.
Define
$$\Gamma:=\phi^{-1}({\R}P^1).$$
This is a real graph on ${\C}P^1$ (invariant with respect
to the complex conjugation) and which contains ${\R}P^1$.
Each vertex of $\Gamma$ has even valency, and the multiplicity of a critical point with real critical value
of $\phi$ is half its valency. The graph $\Gamma$ contains the inverse images of $-1$, $\infty$ and $0$ which are the sets
of roots of $P$, $Q$ and $f$, respectively. Denote by the same letter $p$ (resp. $q$ and $r$) the
points of $\Gamma$ which are mapped to $-1$ (resp. $\infty$ and $0$).
Orient the real axis on the target space via the arrows $-1 \rightarrow \infty \rightarrow 0 \rightarrow -1$
(orientation given by the decreasing order in $\R$) and pull back this orientation by $\phi$.
The graph $\Gamma$ becomes an oriented graph,
with the orientation given by arrows $p \rightarrow q \rightarrow r \rightarrow p$.

Clearly, the arrangement of the real roots of $f$, $P$ and $Q$ together with their multiplicities can be extracted from
the graph $\Gamma$. We encode this arrangement together with the multiplicities
by what is called a root scheme.

\begin{defi}[\cite{Br}]
A root scheme is a $k$-uple $((l_1,m_1),\dotsc,(l_k,m_k)) \in (\{p,q,r\} \times \N)^k$.
A root scheme is realizable by polynomials of degree $d$ if there exist real polynomials
$P$ and $Q$ such that $f=P-Q$ has degree $d$ and if
$\rho_1< \dotsc < \rho_k$ are the real roots of $f$, $P$ and $Q$, then $l_i=p$ (resp. $q$, $r$)
if $\rho_i$ is root of $P$ (resp. $Q$, $f$) and $m_i$ is the multiplicity of $\rho_i$.
\end{defi}

Conversely, suppose we are given a real graph $\Gamma \subset {\C}P^1$
together with a real continuous map $\varphi: \Gamma \rightarrow {\R}P^1$.
Denote the inverse images of $-1$, $\infty$ and $0$ by letters $p$, $q$ and $r$, respectively, and
orient $\Gamma$ with the pull back by $\varphi$ of the above orientation of ${\R}P^1$.
This graph is called {\it a real rational graph}~\cite{Br} if
\begin{itemize}
\item any vertex of $\Gamma$ has even valence,
\item any connected component of ${\C}P^1 \setminus \Gamma$ is homeomorphic to an open disk,
\end{itemize}
Then, for any connected component $D$ of ${\C}P^1 \setminus \Gamma$, the map $\varphi_{| \partial D}$
is a covering of ${\R}P^1$ whose degree $d_D$ is the number of letters $p$ (resp. $q$, $r$)
in $\partial D$. We define the {\it degree} of $\Gamma$ to be half the sum of the degrees $d_D$
over all connected components of ${\C}P^1 \setminus \Gamma$. Since $\varphi$ is a real map,
the degree of $\Gamma$ is also the sum of the degrees $d_D$ 
over all connected components $D$ of ${\C}P^1 \setminus \Gamma$ contained in one connected component of ${\C}P^1 \setminus {\R}P^1$.

\begin{prop}[\cite{Br,O}] \label{P:rootgraph}
A root scheme is realizable by polynomials of degree $d$ if and only if it
can be extracted from a real rational graph of degree $d$ on ${\C}P^1$.
\end{prop}

Let us explain how to prove the if part in Proposition~\ref{P:rootgraph} (see~\cite{Br,Brt,O,NSV}).
For each connected component $D$ of ${\C}P^1 \setminus \Gamma$, extend $\varphi_{|\partial D}$ to a, branched if $d_D>1$,
covering of one connected component of ${\C}P^1 \setminus {\R}P^1$, so that two adjacent connected components of
${\C}P^1 \setminus \Gamma$ project to differents connected components of ${\C}P^1 \setminus {\R}P^1$.
Then, it is possible to glue continuously
these maps in order to obtain a real branched covering $\varphi:{\C}P^1 \rightarrow {\C}P^1$ of degree $d$.
The map $\varphi$ becomes a real rational map of degree $d$ for the standard complex structure
on the target space and its pull-back by $\varphi$ on the source space.
There exist then real polynomials $P$ and $Q$ such that $f=P-Q$ has degree $d$ and
$\varphi =f/Q=P/Q-1$, so that the points $p$ (resp. $q$, $r$)
correspond to the roots of $P$ (resp. $Q$, $f$) and $\Gamma=\varphi^{-1}({\R}P^1)$.

As in~\cite{Br}, we will abbreviate a sequence $S$ repeated $u$ times in a root scheme by $S^u$.
If $u=0$, then $S^u$ is the empty sequence.

\subsection{Constructions}

We are going to prove the existence of root schemes by constructing real rational graphs $\Gamma$ on ${\C}P^1$.
Since these graphs are real, only half of them will be drawn and the horizontal line will represent the
real axis ${\R}P^1$.
\smallskip

\begin{prop}\label{P:rootschemeR=0}

For any even integer $n=2k>0$, the root scheme
$$\left([(q,2 ),(p,2)]^{k}, (q,1), (r,1)^{2k+1}, (p,1) \right)$$
is realizable by polynomials of degree $n+1$.

For any odd integer $n=2k+1>0$, the root scheme
$$\left([(q,2 ),(p,2)]^{k}, (q,2),(p,1)(r,1)^{2k+2},(p,1)\right)$$
is realizable by polynomials of degree $n+1$.
\end{prop}

\noindent {\bf Proof.}
According to Proposition~\ref{P:rootgraph}, it suffices to construct
a real rational graph $\Gamma_{n+1}$ of degree $n+1$ on ${\C}P^1$ from which the desired root scheme can be extracted.
This is done in Figure 1 for $n=2,4$ and in Figure 2 for $n=3,5$. These figures provide the induction step
$n \rightarrow n+2$ for constructing suitable graphs $\Gamma_{n+1}$ for any positive integer $n$.
$\Box$ 
\smallskip

$$
\begin{picture}(0,0)%
\includegraphics{BihRojfig2.pstex}%
\end{picture}%
\setlength{\unitlength}{3750sp}%
\begingroup\makeatletter\ifx\SetFigFont\undefined%
\gdef\SetFigFont#1#2#3#4#5{%
  \reset@font\fontsize{#1}{#2pt}%
  \fontfamily{#3}\fontseries{#4}\fontshape{#5}%
  \selectfont}%
\fi\endgroup%
\begin{picture}(5496,3526)(-197,919)
\end{picture}%

$$
\begin{center}
{\sc Figure 1.}
\end{center}
$$
\begin{picture}(0,0)%
\includegraphics{BihRojfig1.pstex}%
\end{picture}%
\setlength{\unitlength}{3750sp}%
\begingroup\makeatletter\ifx\SetFigFont\undefined%
\gdef\SetFigFont#1#2#3#4#5{%
  \reset@font\fontsize{#1}{#2pt}%
  \fontfamily{#3}\fontseries{#4}\fontshape{#5}%
  \selectfont}%
\fi\endgroup%
\begin{picture}(6359,3732)(-11,1048)
\end{picture}%

$$
\begin{center}
{\sc Figure 2.}
\end{center}

Let $n$ and $R$ be integers such that $0<R \leq n$. If $n-R$ is even, define $I$ and $J$ by
$$
\begin{array}{lll}
I & := & [(p,2),(q,2)]^{\frac{n-R}{2}},(q,3),(r,1),(p,2),(r,1)^{n-R+1}\\
 & & \\
J & := &\left\{\begin{array}{ll}
(r,1), (p,1) & \mbox{if $R=1$} \\
 &  \\
(r,1)^{R},(q,2),[(r,1),(p,4),(r,1),(q,4)]^{\frac{R}{2}-1},(r,1),(p,3)& \mbox{if $R$ is even}\\
 &  \\
(r,1)^{R},(p,2),[(r,1),(q,4),(r,1),(p,4)]^{\frac{R-3}{2}},(q,4),(r,1),(p,3) & \mbox{if $R$ is odd}
\end{array}
\right.
\end{array}$$
If $n-R$ is odd, define
$$
\begin{array}{lll}
I & := & [(p,2),(q,2)]^{\frac{n-R-1}{2}},(p,2),(q,3),(r,1),(p,2),(r,1)^{n-R+1}\\
 & & \\
J & := & \left\{\begin{array}{ll}
(r,1), (q,1) & \mbox{if $R=1$} \\
 &  \\
(r,1)^{R},(p,2),[(r,1),(q,4),(r,1),(p,4)]^{\frac{R}{2}-1},(r,1),(q,3)& \mbox{if $R$ is even}\\
 &  \\
(r,1)^{R},(q,2),[(r,1),(p,4),(r,1),(q,4)]^{\frac{R-3}{2}},(p,4),(r,1),(q,3) & \mbox{if $R$ is odd}
\end{array}
\right.
\end{array}$$
Note that the definitions of $J$ with $n-R$ odd and with $n-R$ even are obtained from each other
by permuting $p$ and $q$.
\smallskip

$$
\begin{picture}(0,0)%
\includegraphics{Bihfig3.pstex}%
\end{picture}%
\setlength{\unitlength}{2565sp}%
\begingroup\makeatletter\ifx\SetFigFont\undefined%
\gdef\SetFigFont#1#2#3#4#5{%
  \reset@font\fontsize{#1}{#2pt}%
  \fontfamily{#3}\fontseries{#4}\fontshape{#5}%
  \selectfont}%
\fi\endgroup%
\begin{picture}(9024,7373)(-1811,-10203)
\end{picture}%

$$
\begin{center}
{\sc Figure 3: graph $\Gamma_I$.}
\end{center}
\smallskip

\begin{prop}\label{P:rootschemeRnot0}
For any integers $n$ and $R$ such that $0<R \leq n$, the root scheme $(I \, , \,  J)$ is realizable by polynomials
of degree $n+R+1$.
\end{prop}

\noindent {\bf Proof.}
First, we contruct graphs $\Gamma_I$ and $\Gamma_j$, which are not real rational graphs, but from which
the sequences $I$ and $J$ can be extracted, respectively.
Figure 3 shows $\Gamma_I$ for $n-R=0,1,2,3$, and indicates how
to construct $\Gamma_I$ for any integer value of $n-R$ with $0 \leq n-R <n$.
Similarly, Figure 4 shows $\Gamma_J$ for $R=1,2,3,4$ and $n-R$ even, and indicates 
how to construct $\Gamma_J$ for any integers $n$ and $R$ such that $0<R \leq n$ and $n-R$ is even.
The graph $\Gamma_J$ for $n-R$ odd is obtained from the graph $\Gamma_J$ with $n-R$ even
and the same value of $R$ by permuting $p$ and $q$ and reversing all the arrows.
For any integers $n$ and $R$ such that $0<R \leq n$
we can glue the corresponding $\Gamma_I$ and $\Gamma_J$ in order
to obtain a real rational graph of degree $n+R+1$ whose root
scheme is  $(I \, , \,  J)$ (see Figure 5 for $n=4$ and $R=3$).
The result follows then from Proposition~\ref{P:rootgraph}.
$\Box$
\smallskip

$$
\begin{picture}(0,0)%
\includegraphics{Bihfig4.pstex}%
\end{picture}%
\setlength{\unitlength}{2565sp}%
\begingroup\makeatletter\ifx\SetFigFont\undefined%
\gdef\SetFigFont#1#2#3#4#5{%
  \reset@font\fontsize{#1}{#2pt}%
  \fontfamily{#3}\fontseries{#4}\fontshape{#5}%
  \selectfont}%
\fi\endgroup%
\begin{picture}(9024,7538)(-911,-10053)
\end{picture}%

$$
\begin{center}
{\sc Figure 4: graph $\Gamma_J$.}
\end{center}
\smallskip

$$
\begin{picture}(0,0)%
\includegraphics{Bihfig6.pstex}%
\end{picture}%
\setlength{\unitlength}{2565sp}%
\begingroup\makeatletter\ifx\SetFigFont\undefined%
\gdef\SetFigFont#1#2#3#4#5{%
  \reset@font\fontsize{#1}{#2pt}%
  \fontfamily{#3}\fontseries{#4}\fontshape{#5}%
  \selectfont}%
\fi\endgroup%
\begin{picture}(11424,1538)(-2711,-5853)
\end{picture}%

$$
\begin{center}
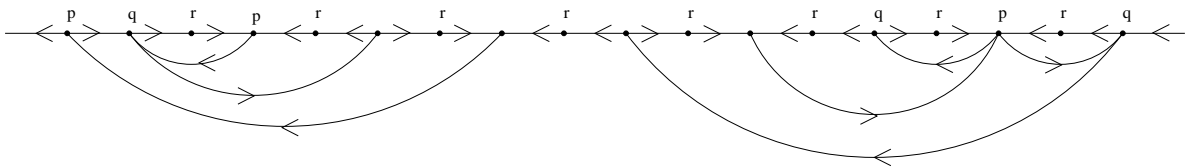
{\sc Figure 5: gluing of $\Gamma_I$ and $\Gamma_J$ for $n=4$, $R=3$.}
\end{center}

\subsection{Sharpness of bounds}
\begin{prop}\label{P:positiveexistence}
For any positive integer $n$, there exist a non degenerate circuit $\calC \subset {\Z}^n$
and a generic real polynomial system supported on $\calC$ whose number of positive solutions
is equal to $n+1$.
\end{prop}

\noindent {\bf Proof.}
We prove only the case $n$ even since the case $n$ odd is similar (the case $n=1$ is trivial).
So assume that $n=2k>0$ is even.
According to Proposition~\ref{P:rootschemeR=0},
there exist polynomials $P$ and $Q$ such that $f=P-Q$ has degree $n+1$
and the corresponding root scheme is  
$$\left([(q,2), (p,2)]^k, (q,1), (r,1)^{2k+1}, (p,1)\right).$$

Composing on the right $\phi=f/Q=P/Q-1$ with a real automorphism of ${\C}P^1$ (a real rational map of degree $1$),
we can choose three values for the points is in this root scheme.
Reading the above scheme from the left to the right, we choose
$p=0$ for the last but one $p$, a positive value for the last $q$, and $p=\infty$
for the last $p$.
This gives the following arrangement of the roots of $P$, $Q$ and $f$
(where we naturally omit the root at infinity)
$$q_1 < p_1 < \dotsc < q_k < p_k=0  < q_{k+1} < r_1 <r_2 < \dotsc < r_{2k+1},$$
so that $f=P-Q$ with
$$
P(x)=a x^2(x-p_1)^2 \cdots (x-p_{k-1})^2 \; , \quad
Q(x)=b (x-q_{k+1}) (x-q_1)^2 \cdots (x-q_k)^2,$$
and $f=P-Q$ has $n+1=2k+1$ positive roots which are bigger than the roots of
$P$ and $Q$. Note that $ab >0$ for otherwise the polynomial $f=P-Q$ could not have roots bigger than the roots of $P$ and $Q$.
Dividing $f$ by $a$ and setting $g_1(x)=x-p_1, \dotsc, g_{k-1}=x-p_{k-1}$, $g_k(x)=x-q_1,\dotsc, g_{n-1}(x)=x-q_k$ and
$g_n(x)=b(x-q_{k+1})/a$, we obtain a polynomial $f$ of the form~\eqref{E:eliminant} with ${\ell}=1$,
the exponents are coprime, and such that ${\calR}_+$ has $n+1$ elements.
Let us now give explicitely a non degenerate circuit $\calC \subset {\Z}^n$
and then a system~\eqref{E:reduced} whose eliminant is $f$.
Take the circuit $\calC$ with $w_{-1}=0$, $w_0=e_1$, $w_1=e_{2},\dotsc,w_{n-1}=e_{n}$
and $w_n=2(e_1+e_2+\dotsc+e_{k}-e_{k+1}-\dotsc-e_{n})$. The primitive affine relation
on $\{e_1,w_1,\dotsc,w_n\}$ is
$$2e_1+2(w_1+\dotsc +w_{k-1})=2(w_k+\dotsc+w_{n-1})+w_n$$
and has the desired coefficients.
Hence $f(x_1)$ is the eliminant of the system $x^{w_i}=g_i(x_1)$, $i=1,\dotsc,n$, given explicitely by
$$
\left\{
\begin{array}{l}
x_2 =  x_1-p_1 \quad , \; \dotsc \; , \; x_{k} = x_1-p_{k-1} \\
x_{k+1} =  x_1-q_1 \; , \; \dotsc \; , \; x_{n} = x_1-q_{k}     \\

{\left(
\frac{x_1 \cdots x_{k}}{x_{k+1} \cdots x_{n}}
\right)}^2  =  b(x_1-q_{k+1})/a.
\end{array}
\right.
$$
According to Proposition~\ref{P:positive}, this system has $n+1$ positive solutions.
$\Box$
\smallskip

It may be interesting to note that the previous system has in fact only these $n+1$ positive solutions
as solutions in the complex torus. To see this, one can compute that the volume of the circuit is $n+1$,
and use Kouchnirenko theorem. 
\smallskip
 
\begin{prop}\label{P:existence}
For any integers $n$ and $R$ such that $0 \leq R \leq n$,
there exist a non degenerate circuit $\calC \subset {\Z}^n$
with $rk \, (\bar{\calC})=R$, and a generic real polynomial system supported on $\calC$ whose number of real solutions
is equal to $2^{n-R}(n+R+1)$.
\end{prop}

\noindent {\bf Proof.}
Note that for $R=0$ this follows from Proposition~\ref{P:positive} (see the proof of Theorem~\ref{T:upper bound}).
We will give the proof only when $n$ and $R$ are even since the other cases are similar (the case
$n=1$ is trivial).

So let $n$ and $R$ be even integers such that $0 < R \leq n$.
According to Proposition~\ref{P:rootschemeRnot0},
there exist polynomials $P$ and $Q$ such that $f=P-Q$ has degree $n+R+1$
and the corresponding root scheme is  
{\small
$$
\left([(p,2),(q,2)]^{\frac{n-R}{2}},(q,3),(r,1),(p,2),(r,1)^{n+1},
(q,2),[(r,1),(p,4),(r,1),(q,4)]^{\frac{R}{2}-1},(r,1),(p,3) \right).
$$
}
We may choose three values for the roots in this root scheme.
Choose $p=0$ for the $p$ in the sequence $(q,3),(r,1),(p,2)$, a negative value for the $r$ in the same sequence,
and the infinity value for the last $p$ on the right in the root scheme.
Then, the polynomials $P$ and $Q$ can be written as
$$
P(x)=x^2 \prod_{i=1}^{\frac{n-R}{2}} (g_i(x))^2 \cdot \prod_{i=\frac{n-R}{2}+1}^{\frac{n}{2}-1}((g_i(x))^4,
$$
$$Q(x)=(g_{\frac{n}{2}}(x))^3 \cdot \prod_{i=\frac{n}{2}+1}^{n-\frac{R}{2}+1} (g_i(x))^2 \cdot
\prod_{i=n-\frac{R}{2}+2}^{n}((g_i(x))^4
$$
so that 
\medskip

\begin {tabular}{ll}
($\star$)$\;$ & {\it $f=P-Q$ has $n+R+1$ (non zero) real roots, all but one are positive, and}\\
  & {\it $g_i$ is positive at the real roots of $f$
for $i=1,\dotsc,\frac{n-R}{2}$ and $i=\frac{n}{2},\dotsc,n-\frac{R}{2}$.}
\end{tabular}
\medskip

\noindent Here the $g_i$ have all, but one, the form, $g_i(x)=x-\rho_i$ with $\rho_i$ a finite non zero root
$p$ if $i < \frac{n}{2}$ and a root
$q$ otherwise. The last one is of the form $c(x-\rho_i)$ where $c$ is a positive real number due to the fact
$f$ has a real root bigger than the roots of $P$ and $Q$ (see the proof of Proposition~\ref{P:positiveexistence}).
The polynomials $g_i$ mentioned in ($\star$) correspond to  the sequence
$[(p,2),(q,2)]^{\frac{n-R}{2}},(q,3)$ on the left of the above root scheme.
\smallskip

Now, define $v_{n-\frac{R}{2}+1}=-(e_2+\dotsc+e_{n-1})-3e_n$, $v_{\frac{n}{2}}=2(e_2+\dotsc+e_n)$, and
$$
\begin{array}{lll}
v_i & = & 2e_{i+1} \;,  \quad i=1,\dotsc,\frac{n-R}{2}\\
 & & \\
v_i & = & e_{i+1} \;,  \quad i=\frac{n-R}{2}+1,\dotsc,\frac{n}{2}-1\\
 & & \\
v_i & = & -2e_{i} \;,  \quad i=\frac{n}{2}+1,\dotsc,n-\frac{R}{2}\\
 & & \\
v_i & = & -e_{i-1} \;,   \quad i=n-\frac{R}{2}+2,\dotsc,n 
\end{array}
$$
The primitive affine relation on $\{0,v_1,\dotsc,v_n\}$ is the desired one
$$2\sum_{i=1}^{\frac{n-R}{2}}v_i+4\sum_{i=\frac{n-R}{2}+1}^{\frac{n}{2}-1}v_i=
3v_{\frac{n}{2}}+2\sum_{i=\frac{n}{2}+1}^{n-\frac{R}{2}+1}v_i+4\sum_{i=n-\frac{R}{2}+2}^{n}v_i,
$$
and it's easy to see that the rank of the matrix $\bar{B}=(\bar{v}_1, \dotsc, \bar{v}_n)$
is $R-1$. We want to find integers
$l_1,\dotsc,l_n$ verifying the relation~\eqref{E:primitiveaffinerelationl}
and such that the sign conditions in Proposition~\ref{P:real} are verified at each root $r$ of $f$.
Here, these sign conditions are given by the $v_i$ which belong to $2{\Z}^n$
\begin{equation}\label{E:signproof}
\begin{array}{lll}
g_i(r)/r^{l_i} &  >  0& ,  \quad i=1,\dotsc,\frac{n-R}{2} \quad \mbox{and} \quad 
 i=\frac{n}{2},\dotsc,n-\frac{R}{2}
\end{array}
\end{equation}
In view of property $(\star)$, we want to have $l_i$ even for $i=1,\dotsc,\frac{n-R}{2}$
and $i=\frac{n}{2},\dotsc,n-\frac{R}{2}$ since one root of $f$ is negative.
In fact we can find integers $l_i$ which are all even. Indeed, in our situation,
the relation~\eqref{E:primitiveaffinerelationl} is

\begin{equation}\label{E:relationproof}
2+2\sum_{i=1}^{\frac{n-R}{2}}l_i+4\sum_{i=\frac{n-R}{2}+1}^{\frac{n}{2}-1}l_i=
3l_{\frac{n}{2}}+2\sum_{i=\frac{n}{2}+1}^{n-\frac{R}{2}+1}l_i+4\sum_{i=n-\frac{R}{2}+2}^{n}l_i.
\end{equation}
Replacing the first $2$ on the left by $1$, the existence of integers $l_i$ verifying the resulting relation
is due to the fact that $2,4,3$ are coprime. Multiplying then by $2$ the members of this modified
relation gives then the existence of even integers $l_i$ verifying~\eqref{E:relationproof}.

So, let us choose even integers $l_i$ verifying~\eqref{E:relationproof} so that the sign
conditions~\eqref{E:signproof} are satisfied. The proof is almost finished.

Consider the system
\begin{equation}\label{E:systemproof}
x^{w_i}=g_i(x_1) \; , \quad i=1,\dotsc,n,
\end{equation}
where $w_i:=l_ie_1+v_i$.
This system is equivalent to a (generic) system supported on the circuit $\calC=\{w_{-1}=0,w_0=e_1,w_1,\dotsc,w_n\}$.
The primitive affine relation on $\calC$ is
$$
2w_0+2\sum_{i=1}^{\frac{n-R}{2}}w_i+4\sum_{i=\frac{n-R}{2}+1}^{\frac{n}{2}-1}w_i=
3w_{\frac{n}{2}}+2\sum_{i=\frac{n}{2}+1}^{n-\frac{R}{2}+1}w_i+4\sum_{i=n-\frac{R}{2}+2}^{n}w_i,
$$
so that $f$ is the eliminant of~\eqref{E:systemproof}. The rank modulo $2$ of $\calC$ is
${{rk \, ({\bar \calC})}}=rk \,(\bar{B})+1=R$. The eliminant $f$ has $n+R+1$ real roots and the sign
conditions~\eqref{E:signproof} are satisfied at each real root of $f$.
Therefore, according to Proposition~\ref{P:real} , the system~\eqref{E:systemproof} has $2^{n-R}(n+R+1)$ real solutions.
$\Box$
\smallskip

The case $R=n$ in Proposition~\ref{P:existence} has already been obtained in~\cite{BeBiS}
by different methods.

\begin{thm}\label{T:positivesharp}
For any integers $n, m$ such that $1 \leq m \leq n$, there exist a circuit $\calC \subset {\Z}^n$
with $m(\calC)=m$
and a generic real polynomial system supported on $\calC$ whose number of positive solutions
is $m+1$
\end{thm}

\noindent {\bf Proof.}
By Proposition~\ref{P:positiveexistence}, there exists a circuit
${\calC}'=\{0,w_0=\ell e_1,w_1,\dotsc,w_{m}\} \subset {\Z}^{m}$
with $m({\calC}')=m$ and polynomials $g_1,\dotsc,g_{m}$ such that the system
$x^{w_i}=g_i(x_1)$, $i=1,\dotsc, m$, has ${m}+1$ positive solutions.
By Proposition~\ref{P:positive}, this means that the eliminant $f$ of this system
has $m+1$ positive roots at which $g_1,\dotsc,g_m$ are positive.
Define
$$w_i:=e_i \; , \quad i={m} +1,\dotsc,n,$$
and choose polynomials $g_{m+1},\dotsc,g_n$ (as in~\eqref{E:reduced})
which are positive at the real roots of $f$.
Then, the system $x^{w_i}=g_i(x_1)$, $i=1,\dotsc, n$,
is equivalent to a system supported on the circuit
${\calC}=\{0,w_0=\ell e_1,w_1,\dotsc,w_n\} \subset {\Z}^{n}$. The eliminant of this system is also $f$
and we have $m({\calC})=m({{\calC}'})=m$. 
Proposition~\ref{P:positive} and our choice of $g_{{m}+1},\dotsc,g_n$ implies
that this system has also ${m}+1$ positive solutions.
$\Box$
\smallskip

\begin{thm}\label{T:sharp}
For any integers $n, m, R$ such that $0 \leq R \leq n$
and $1 \leq m \leq n$, there exist a circuit $\calC \subset {\Z}^n$
with $m(\calC)=m$, ${rk \, ({\bar \calC})}=R$,
and a generic real polynomial system supported on $\calC$ whose number of real solutions,
$N$, verifies:

\begin{enumerate}
\item if $R \leq m$, then
$$N= 2^{n-{R}} \cdot \left({m}+ {R}+1\right).$$
\item if ${R} \geq {m}$, then
$$N = 2^{n-{R}} \cdot \left(2{m}+1\right).$$
\end{enumerate}
\end{thm}

\noindent {\bf Proof.}
The case $R=0$ is an immediate consequence of Theorem~\ref{T:positivesharp}.
Suppose that $0 < R \leq {m}$.
By Proposition~\ref{P:existence}, there exist a circuit
${\calC}'=\{0,w_0=\ell e_1,w_1,\dotsc,w_{m}\} \subset {\Z}^{m}$
with $m({\calC}')=m$, $rk \, (\bar{{\calC}'})=R$,
and polynomials $g_1,\dotsc,g_{m}$ such that the system
\begin{equation}\label{E:system1}
x^{w_i}=g_i(x_1) \; , \quad i=1,\dotsc, m,
\end{equation}
has $2^{m-{R}} \cdot \left({m}+ {R}+1\right)$
real solutions. Here and in the rest of the proof $g_i$ is a polynomial as in~\eqref{E:reduced}.
Since $ rk \, (\bar{{\calC}'}) \neq 0$, we may assume that $\ell$ is odd.
Define
$$w_i:=2e_i \; , \quad i={m} +1,\dotsc,n,$$
and choose polynomials $g_{m+1},\dotsc,g_n$
which are positive at the real roots of the eliminant $f$ of~\eqref{E:system1}.
Then, the system
\begin{equation}\label{E:system2}
x^{w_i}=g_i(x_1) \; , \quad i=1,\dotsc,n,
\end{equation}
is equivalent to a system supported on the circuit
${\calC}=\{0,w_0=\ell e_1,w_1,\dotsc,w_n\} \subset {\Z}^{n}$.
We have 
$m(\calC)=m({\calC'})=m$, ${rk \, ({\bar {\calC}})}=rk \, (\bar{{\calC}'})=R$
and $f$ is also the eliminant of~\eqref{E:system2}.
Moreover the sign conditions in Proposition~\ref{P:real} which correspond to
~\eqref{E:system2} are obtained from those corresponding to~\eqref{E:system1}
by adding $g_i(r)>0$ for $i=m+1,\dotsc,n$ (note that $l_{m+1}=\dotsc=l_n=0$).
But since~\eqref{E:system1} has $2^{m-{R}} \cdot \left({m}+ {R}+1\right)$ real solutions,
and $\ell$ is odd, it follows from Proposition~\ref{P:real} that the sign conditions
corresponding to~\eqref{E:system1} are satisfied at ${m}+ {R}+1$ real roots of $f$.
By our choice of $g_{m+1},\dotsc,g_n$, the sign conditions corresponding
to~\eqref{E:system2} are satisfied at the same ${m}+ {R}+1$ real roots of $f$.
Hence, again by Proposition~\ref{P:real}, the system~\eqref{E:system2} has
$2^{n-{R}} \cdot \left({m}+ {R}+1\right)$ real solutions.

Finally, suppose that ${R}> {m}$. The idea is exactly as before.
By Proposition~\ref{P:existence}, there exist a circuit
${\calC}'=\{0,w_0=\ell e_1,w_1,\dotsc,w_{m}\} \subset {\Z}^{m}$
with $m({\calC}')=m$, $rk \, (\bar{{\calC}'})=m$ and polynomials $g_1,\dotsc,g_{m}$ such that the system~\eqref{E:system1}
has $2{m}+1$ real solutions. According to Proposition~\ref{P:real} (see also Remark~\ref{R:redundant}),
this means that the eliminant $f$ of this system has $2{m}+1$ real roots.
Define
$$w_i:=e_i \; , \quad i={m} +1,\dotsc,{R}, \quad w_i:=2e_i \; , \quad i={R} +1,\dotsc,n$$
and choose polynomials $g_{R+1},\dotsc,g_n$
which are positive at the roots of the eliminant $f$.
Then, the system~\eqref{E:system2}
is equivalent to a system supported on the circuit
${\calC}=\{0,w_0=\ell e_1,w_1,\dotsc,w_n\} \subset {\Z}^{n}$.
We have $m(\calC)={m}$,
$ rk \, (\bar{\calC})=R$ and~\eqref{E:system2} has
$2^{n-R} \cdot \left(2{m}+1\right)$ real solutions.
$\Box$
\smallskip


\providecommand{\bysame}{\leavevmode\hbox to3em{\hrulefill}\thinspace}
\providecommand{\MR}{\relax\ifhmode\unskip\space\fi MR }
\providecommand{\MRhref}[2]{%
  \href{http://www.ams.org/mathscinet-getitem?mr=#1}{#2}
}
\providecommand{\href}[2]{#2}

\end{document}